# Integral representations of the Legendre chi function

## Djurdje Cvijović[1]

[1] Atomic Physics Laboratory, Vinča Institute of Nuclear Sciences,

P.O. Box 522, 11001 Belgrade, Serbia

E-mail: **djurdje@vinca.rs**

**Abstract.** We, by making use of elementary arguments, deduce integral representations of the Legendre chi function $\chi_s(z)$ valid for $|z| < 1$ and $\operatorname{Re} s > 1$. Our earlier established results on the integral representations for the Riemann zeta function $\zeta(2n+1)$ and the Dirichlet beta function $\beta(2n)$, $n \in N$, are direct consequence of these representations.

## 1. Introduction

In what follows, we shall refer to the following series

$$\chi_s(z) = \sum_{k=0}^{\infty} \frac{z^{2k+1}}{(2k+1)^s} \qquad (|z| \leq 1; \operatorname{Re} s > 1). \qquad (1)$$

as the Legendre chi function $\chi_s(z)$. In the particular case when $s = n = 2, 3, \ldots$, it is usually called Legendre's chi function of order n. (for details see [1, pp. 17-19] and [2, pp. 282-283])

$\chi_s(z)$ resembles the Dirichlet series for the polylogarithm function $Li_s(z)$ [3, Appendix II. 5, pp. 762-763],

$$Li_s(z) = \sum_{k=1}^{\infty} \frac{z^k}{k^s} \qquad (|z| < 1; s \in C)$$

and, indeed, is trivially expressible in terms of the polylogarithm as

$$\chi_s(z) = \frac{1}{2}\left(L_s(z) - L_s(-z)\right).$$

Adrien-Marie Legendre was first to study the chi function in his famous book "Exercices de calcul integral" (1811) and he actually used $\phi$ to denote it. The present notation goes back to Edward's classical text "A Treatise on the Integral Calculus" (1889). A good review of the theory of the Legendre chi function $\chi_n(z)$ of order n = 2, 3, 4, … is given by Lewin in his two standard texts [1] and [2]. Berndt's treatise can serve as an excellent introductory text on this (and numerous related functions) and as an encyclopaedic source [4]. For newer results concerning $\chi_s(z)$ we refer to Boersma and Dempsey [5] and Cvijović and Klinowski [6].

Here, by making use of elementary arguments we deduce new definite integral definitions of the Legendre chi function. Our earlier established results [7, Theorem 1] on the integral representations for the Riemann zeta function $\zeta(2n+1)$ and the Dirichlet beta function $\beta(2n)$, $n \in N$, directly follow as corollaries of these representations.

## 2. Statement of the results

First, we need to introduce several definitions. The Riemann zeta function $\zeta(s)$ is defined by [8, p. 807, Eqs. 23.2.1]:

$$\zeta(s) = \sum_{k=1}^{\infty} \frac{1}{k^s} \qquad (\operatorname{Re} s > 1) \qquad (2)$$

and can be continued meromorphically to the whole *s*-plane with a simple pole at *s* = 1.

The Dirichlet lambda and beta function, $\lambda(s)$ and $\beta(s)$, are for $\operatorname{Re} s > 1$ defined as (compare with [8, p. 807, Eq. 23.2.20])

$$\lambda(s) = \sum_{k=0}^{\infty} \frac{1}{(2k+1)^s} = (1 - 2^{-s})\zeta(s) \qquad (3)$$

and

$$\beta(s) = \sum_{k=0}^{\infty} \frac{(-1)^k}{(2k+1)^s}, \tag{4}$$

(compare with [8, p. 807, Eq. 23.2.21]) respectively. It should be noted that there exists relation between the Riemann zeta and the Dirichlet lambda function.

It is convenient to consider the functions defined by means of the following series

$$S(s,x) = \sum_{k=0}^{\infty} \frac{\sin(2k+1)x}{(2k+1)^s} \tag{5}$$

and

$$C(s,x) = \sum_{k=0}^{\infty} \frac{\cos(2k+1)x}{(2k+1)^s}, \tag{6}$$

where $\operatorname{Re} s > 1$. We remark that the definitions of $S(s,x)$ and $C(s,x)$ ensure the convergence of each of the series involved: both series converge uniformly for all real values of x when $\operatorname{Re} s > 1$.

We also use the Euler polynomials of degree n in x, $E_n(x)$, defined by means of the following generating function [8, p. 804, Eq. 23.1.1]

$$\frac{2te^{tx}}{e^t + 1} = \sum_{n=0}^{\infty} E_n(x) \frac{t^n}{n!} \qquad (|t| < \pi). \tag{7}$$

Our results are as follows.

**Theorem 1.** Set δ = 1 and δ = ½. Assume that s and z are complex numbers, let $\chi_s(z)$ be the Legendre chi function and let $S(s,x)$ and $C(s,x)$ be defined as in Eqs. (5) and (6), respectively. If $\operatorname{Re} s > 1$ and $|z| < 1$, then:

$$\chi_s(z) = \frac{1}{\delta}\int_0^\delta S(s,\pi t)\frac{2z(1+z^2)\sin(\pi t)}{1-2z^2\cos(2\pi t)+z^4}\,dt, \qquad (8a)$$

and

$$\chi_s(z) = \frac{1}{\delta}\int_0^\delta C(s,\pi t)\frac{2z(1-z^2)\cos(\pi t)}{1-2z^2\cos(2\pi t)+z^4}\,dt. \qquad (8b)$$

**Theorem 2.** Assume that n is a positive integer, let $\chi_n(z)$ be the Legendre chi function and let $E_n(x)$ be the Euler polynomials of degree n in x given by Eq. (7). Then:

$$\chi_{2n+1}(z) = (-1)^n \frac{1}{\delta}\frac{\pi^{2n+1}}{4(2n)!}\int_0^\delta E_{2n}(t)\frac{2z(1+z^2)\sin(\pi t)}{1-2z^2\cos(2\pi t)+z^4}\,dt \qquad (9a)$$

$$\chi_{2n}(z) = (-1)^n \frac{1}{\delta}\frac{\pi^{2n}}{4(2n-1)!}\int_0^\delta E_{2n-1}(t)\frac{2z(1-z^2)\cos(\pi t)}{1-2z^2\cos(2\pi t)+z^4}\,dt, \qquad (9b)$$

where $|z| < 1$.

**Corollary.** Assume that n is a positive integer. Let $\lambda(s)$ and $\beta(s)$ be the Dirichlet lambda and beta function defined as in (3) and (4) respectively and let $E_n(x)$ be the Euler polynomials of degree n in x given by Eq. (7). We then have:

$$\lambda(2n+1) = (-1)^n \frac{1}{\delta} \frac{\pi^{2n+1}}{4(2n)!} \int_0^\delta E_{2n}(t)\csc(\pi t)dt, \qquad (10a)$$

$$\beta(2n) = (-1)^n \frac{1}{\delta} \frac{\pi^{2n}}{4(2n-1)!} \int_0^\delta E_{2n-1}(t)\sec(\pi t)dt. \qquad (10b)$$

**Remark 1.** Note that the existence of the integrals in (10) is assured since the integrands on [0, a], 0 < a < 1, have only removable singularities. This can be easily demonstrated by making use of some basic properties of $E_n(x)$.

**Remark 2.** In view of the relationship between the Riemann zeta and the Dirichlet lambda function (3) it is clear that the expression in (10a) could be rewritten as

$$\zeta(2n+1) = (-1)^n \frac{2^{2n+1}}{(2^{2n+1}-1)} \frac{1}{\delta} \frac{\pi^{2n+1}}{4(2n)!} \int_0^\delta E_{2n}(t)\csc(\pi t)dt,$$

thus giving the integral representation of the Riemann zeta function for odd-integers arguments.

**Remark 3.** We have failed to find in the literature our formulae given by Theorem 1 and Theorem 2. The integral for $\beta(2n)$ in (10b) when δ = 1 is well known [8, p. 807, Eq. 23.2.17] while the other integrals in (10) were recently deduced [7, Theorem 1, p. 436].

## 3. Proof of the results

First of all, we need the following two elementary lemmas.

**Lemma 1.** Suppose that n is a positive integer and that $\alpha = 1$ and $\alpha = \tfrac{1}{2}$. Then the formulae hold:

$$\int_0^{2\alpha\pi} \sin(nt) \frac{2z\sin t}{1-2z\cos t + z^2} dt = 2\alpha\pi z^n, \tag{11a}$$

$$\int_0^{2\alpha\pi} \cos(nt) \frac{1-z^2}{1-2z\cos t + z^2} dt = 2\alpha\pi z^n. \tag{11b}$$

**Proof.** We first establish the case $\alpha = 1$ in (11). Consider the integrals $I_1$ and $I_2$ given respectively by

$$I_1 = \int_0^{2\pi} e^{int} \frac{2z\sin t}{1-2z\cos t + z^2} dt =$$

$$= \int_0^{2\pi} \cos(nt) \frac{2z\sin t}{1-2z\cos t + z^2} dt + i\int_0^{2\pi} \sin(nt) \frac{2z\sin t}{1-2z\cos t + z^2} dt \tag{12a}$$

and by

$$I_2 = \int_0^{2\pi} e^{int} \frac{1-z^2}{1-2z\cos t + z^2} dt =$$

$$= \int_0^{2\pi} \cos(nt) \frac{1-z^2}{1-2z\cos t + z^2} dt + i \int_0^{2\pi} \sin(nt) \frac{1-z^2}{1-2z\cos t + z^2} dt. \qquad (12b)$$

It is not difficult to evaluate these two integrals. To this end we use the contour integration method in the complex plane and calculus of residues. By setting

$$\tau = e^{it}, \quad \cos t = \frac{1}{2}(\tau + 1/\tau) \text{ and } \sin t = \frac{1}{2i}(\tau - 1/\tau),$$

where $i^2 = -1$, we arrive at

$$I_1 = \oint_{|\tau|=1} \frac{i\tau^n(\tau^2-1)}{(\tau-z)(\tau-1/z)} \frac{d\tau}{i\tau} = \oint_{|\tau|=1} \frac{\tau^{n-1}(\tau^2-1)}{(\tau-z)(\tau-1/z)} d\tau =$$

$$= 2\pi i \operatorname*{Res}_{\tau=z} \frac{\tau^{n-1}(\tau^2-1)}{(\tau-z)(\tau-1/z)} = 2\pi i\, z^n, \qquad (13a)$$

and

$$I_2 = \oint_{|\tau|=1} \frac{\tau^{n+1}(z^2-1)}{z(\tau-z)(\tau-1/z)} \frac{d\tau}{i\tau} = \oint_{|\tau|=1} \frac{\tau^n(z^2-1)}{iz(\tau-z)(\tau-1/z)} d\tau =$$

$$= 2\pi i \operatorname*{Res}_{\tau=z} \frac{\tau^n(z^2-1)}{iz(\tau-z)(\tau-1/z)} = 2\pi i(-iz^n) = 2\pi z^n. \qquad (13b)$$

Combining (12) and (13) and equating the real and imaginary parts on both sides gives the integrals in (11) where $\alpha = 1$. Observe that in both integrals in (13) the contour is unit circle and is traversed in the positive (counterclockwise) direction and the only singularities of the integrands that lie inside the contour are at τ = z.

Next, regarding the case α = 1/2 in (10) it follows easily that it reduces to the above considered case α = 1. Indeed, in light of the following well-known property [9, pp. 272-273, Eqs. 2.1.2.20 and 2.1.2.21]

$$\int_0^{2a} f(t)dt = \begin{cases} 0 & f(2a-t) = -f(t) \\ 2\int_0^a f(t)dt & f(2a-t) = f(t) \end{cases} \quad (14)$$

it suffices to verify that both integrands in (11) are such that

$$f(2\pi - t) = f(t).$$

This completes the proof of our lemma.

**Lemma 2.** Assume that k is a nonnegative integer and that δ = 1 and δ = ½. Then we have:

$$\int_0^{\delta\pi} \sin(2k+1)t \frac{2z(1+z^2)\sin t}{1-2z^2\cos(2t)+z^4} dt = \delta\pi z^{2k+1}, \quad (15a)$$

$$\int_0^{\delta\pi} \sin(2k+1)t \frac{2z(1-z^2)\cos t}{1-2z^2\cos(2t)+z^4} dt = \delta\pi z^{2k+1}. \quad (15b)$$

**Proof.** We use Lemma 1 and from (11a) we obtain

$$2\alpha\pi\left(z^n - (-z)^n\right) = \int_0^{2\alpha\pi} \sin(nt) \frac{4z(1+z^2)\sin t}{1-2z^2\cos(2t)+z^4} dt,$$

which, upon setting n = 2 k + 1, k = 0, 1, 2, ..., yields,

$$2\alpha\pi z^{2k+1} = \int_0^{2\alpha\pi} \sin(2k+1)t \frac{2z(1+z^2)\sin t}{1-2z^2\cos(2t)+z^4} dt \qquad (16)$$

where α = 1 and α = ½. However, this formula also holds true for α = ¼, and that can be easily shown by appealing to the property in (14). The required formula (15a) follows at once by substituting α = δ/2 in (16).

Starting from (11b) the formula in (15b) is derived in precisely the same way.

**Proof of Theorem 1.** We use Lemma 2 and from (15a) we have

$$z^{2k+1} = \frac{1}{\delta}\int_0^{\delta} \sin(2k+1)\pi t \frac{2z(1+z^2)\sin(\pi t)}{1-2z^2\cos(2\pi t)+z^4} dt \qquad \delta = 1, 1/2; k = 1, 2, ...,$$

so that

$$\sum_{k=0}^{\infty} \frac{z^{2k+1}}{(2k+1)^s} = \frac{1}{\delta}\int_0^{\delta} \sum_{k=0}^{\infty} \frac{\sin(2k+1)\pi t}{(2k+1)^s} \frac{2z(1+z^2)\sin(\pi t)}{1-2z^2\cos(2\pi t)+z^4} dt \qquad \operatorname{Re} s > 1. (17)$$

Now, we readily find the proposed integral formula in (8a) from the last expression by appealing directly to the definitions of $\chi_s(z)$ in (1) and $S(s,x)$ in (5). Observe that inverting the order of summation and integration in (17) is justified by absolute convergence.

Starting from (15b) we obtain the formula in (8b) in the same manner. This completes our proof.

**Proof of Theorem 2.** First, recall the following Fourier expansions for the Euler polynomials $E_n(x)$ [8, p. 805, Eqs. 23.1.17]

$$E_{2n-1}(x) = (-1)^n \frac{4(2n-1)!}{(2\pi)^{2n}} \sum_{k=1}^{\infty} \frac{\cos(2k+1)\pi x}{(2k+1)^{2n}}$$

where $0 \leq x \leq 1$ for n = 1, 2, …, and [8, p. 805, Eqs. 23.1.18]

$$E_{2n}(x) = (-1)^n \frac{4(2n)!}{(2\pi)^{2n+1}} \sum_{k=1}^{\infty} \frac{\sin(2k+1)\pi x}{(2k+1)^{2n+1}}$$

where $0 \leq x \leq 1$ for n = 1, 2, … .and $0 < x < 1$ for n = 0.

Next, in view of the definitions (5) and (6) the above expansions can be rewritten as

$$S(2n+1, \pi x) = (-1)^n \frac{\pi^{2n+1}}{4(2n)!} E_{2n}(x), \qquad (18a)$$

where $0 \leq x \leq 1$ for n = 1, 2, … .and $0 < x < 1$ for n = 0, and

$$C(2n, \pi x) = (-1)^n \frac{\pi^{2n}}{4(2n-1)!} E_{2n-1}(x), \qquad (18b)$$

where $0 \leq x \leq 1$ for n = 1, 2, …, .

Finally, by employing Theorem 1 in conjunction with (18) we achieve the proposed formulae (9).

**Proof of Corollary.** First, observe that from the definition of $\chi_s(z)$ in (1) we have

$$\chi_s(1) = \sum_{k=0}^{\infty} \frac{1}{(2k+1)^s} = \lambda(s) = (1-2^{-s})\zeta(s) \quad (\text{Re}\, s > 1), \tag{19a}$$

$$\chi_s(i) = \sum_{k=0}^{\infty} \frac{i^{2k+1}}{(2k+1)^s} = i\sum_{k=0}^{\infty} \frac{(-1)^k}{(2k+1)^s} = i\beta(s) \quad (\text{Re}\, s > 1), \tag{19b}$$

where $\zeta(s)$, $\lambda(s)$ and $\beta(s)$ are functions given respectively by (2), (3) and (4) and the symbol $i$ stands for the imaginary unit.

Secondly, it is not difficult to show that

$$\lim_{z \to 1} \frac{2z(1+z^2)\sin(\pi t)}{1-2z^2\cos(2\pi t)+z^4} = \csc(\pi t) \tag{20a}$$

and

$$\lim_{z \to i} \frac{2z(1-z^2)\cos(\pi t)}{1-2z^2\cos(2\pi t)+z^4} = i\sec(\pi t). \tag{20b}$$

Finally, we deduce the integral formulae in (10) starting from (9) and by employing the expressions (19) and (20).